\documentclass{amsart}
\usepackage{latexsym,amsmath,amssymb}
\theoremstyle{definition}
\theoremstyle{remark}

\numberwithin{equation}{section}
\begin{document}

\title[Weighted conditional expectation type operators ]
{Partial normality classes of weighted conditional type operators
on $L^{2}(\Sigma)$ }

\author{\sc\bf Y. Estaremi  }
\address{\sc Y. Estaremi  }
\email{yestaremi@pnu.ac.ir}

\address{Department of Mathematics, Payame Noor University, p. o. box: 19395-3697, Tehran, Iran.}

\thanks{}

\thanks{}



\date{}

\dedicatory{}

\commby{}

%%% ----------------------------------------------------------------------
\begin{abstract}
In this paper, some various partial normality classes of weighted
conditional expectation type operators on $L^{2}(\Sigma)$ are
investigated. Also, some applications of weak hyponormal weighted
conditional type operators are presented.

\noindent {}
\end{abstract}

\maketitle

\section{ \sc\bf Introduction }
The notion of conditional expectation is rightfully thought of as
belonging to the theory of probability. In such a text, it is set
against a background of a probability space $(\Omega,\mathcal{F},
P)$ and $\sigma$-subalgebra ($\sigma$-field as it is commonly
called in probability texts) $\mathcal{G}$ of $\mathcal{F}$. If
$X$ denotes an integrable random variable, then the conditional
expected value of $X$ given $\mathcal{G}$ is the random variable
$E[X|\mathcal{G}]$ such that \\
1. $E[X|\mathcal{G}]$ is $\mathcal{G}$-measurable,\\
2. $E[X|\mathcal{G}]$ satisfies the functional relation
$$\int_{G}E[X|\mathcal{G}]dP=\int_{G}XdP, \ \ \ \ \forall G\in \mathcal{G}.$$
Any number of standard texts will illustrate concisely the
probabilistic formulation and interpretation of the function
$E[X|\mathcal{G}]$, and the reader is invited to consult such
references as \cite{Bill}. Our main interests, however, reside in
the view of conditional expectation as an operator between the
$L^p$-spaces, specially between $L^2$-spaces.\\
Among the earlier investigations along these lines is that of
Shu-Teh Chen Moy in his seminal 1954 paper \cite{mo}. Set within
the familiar framework of a probability space $(\Omega,
\mathcal{F}, P)$, Moy obtains necessary and sufficient conditions
for a linear transformation $T$ between function spaces to be of
the form $TX=E[gX|\mathcal{G}]$, where
$\mathcal{G}\subset\mathcal{F}$ is a $\sigma$-subalgebra and $g$
is a nonnegative measurable function with bounded conditional
expected value. The function $E[gX|\mathcal{G}]$ can best be
described as the weighted conditional expected value of $X$.
Moreover, Conditional expectations have been studied in an
operator theoretic setting, by, for example, R. G. Douglas,
\cite{dou}, de Pagter and Grobler \cite{gd}, P.G. Dodds, C.B.
Huijsmans and B. De Pagter, \cite{dhd}, J. Herron,  \cite{her},
Alan Lambert \cite{lam} and Rao \cite{rao1, rao2}, as positive
operators acting on $L^p$-spaces or Banach function spaces. The
combination of conditional expectation and multiplication
operators appears more often in the service of the study of other
operators rather than being the object of study in and of
themselves.\\ In \cite{e,ej} we investigated some classic
properties of multiplication conditional expectation operators
$M_wEM_u$ on $L^p$ spaces. In this paper we will be concerned with
characterizing weighted conditional expectation type operators on
$L^2(\Sigma)$ in terms of membership in the various partial
normality classes and some applications of them.

\section{ \sc\bf Preliminaries}

Let $(X,\Sigma,\mu)$ be a complete $\sigma$-finite measure space.
For any $\sigma$-subalgebra $\mathcal{A}\subseteq
 \Sigma$, the $L^2$-space
$L^2(X,\mathcal{A},\mu_{\mid_{\mathcal{A}}})$ is abbreviated  by
$L^2(\mathcal{A})$, and its norm is denoted by $\|.\|_2$. All
comparisons between two functions or two sets are to be
interpreted as holding up to a $\mu$-null set. The support of a
measurable function $f$ is defined as $S(f)=\{x\in X; f(x)\neq
0\}$. We denote the vector space of all equivalence classes of
almost everywhere finite valued measurable functions on $X$ by
$L^0(\Sigma)$.

\vspace*{0.3cm} For a $\sigma$-subalgebra
$\mathcal{A}\subseteq\Sigma$, the conditional expectation operator
associated with $\mathcal{A}$ is the mapping $f\rightarrow
E^{\mathcal{A}}f$, defined for all non-negative, measurable
function $f$ as well as for all $f\in L^2(\Sigma)$, where
$E^{\mathcal{A}}f$, by the Radon-Nikodym theorem, is the unique
$\mathcal{A}$-measurable function satisfying
$$\int_{A}fd\mu=\int_{A}E^{\mathcal{A}}fd\mu, \ \ \ \forall A\in \mathcal{A} .$$
As an operator on $L^{2}({\Sigma})$, $E^{\mathcal{A}}$ is
idempotent and $E^{\mathcal{A}}(L^2(\Sigma))=L^2(\mathcal{A})$. If
there is no possibility of confusion, we write $E(f)$ in place of
$E^{\mathcal{A}}(f)$. This operator will play a major role in our
work and we list here some of its useful properties:

\vspace*{0.2cm} \noindent $\bullet$ \  If $g$ is
$\mathcal{A}$-measurable, then $E(fg)=E(f)g$.

\noindent $\bullet$ \ $|E(f)|^2\leq E(|f|^2)$.

\noindent $\bullet$ \ If $f\geq 0$, then $E(f)\geq 0$; if $f>0$,
then $E(f)>0$.

\noindent $\bullet$ \ $|E(fg)|\leq
(E(|f|^2))^{\frac{1}{2}}(E(|g|^{2}))^{\frac{1}{2}}$, (H\"{o}lder
inequality).

\noindent $\bullet$ \ For each $f\geq 0$, $S(f)\subseteq S(E(f))$.

\vspace*{0.2cm}\noindent A detailed discussion and verification of
most of these properties may be found in \cite{rao}.

 Let $f\in L^0(\Sigma)$, then
$f$ is said to be conditionable with respect to $E$ if
$f\in\mathcal{D}(E):=\{g\in L^0(\Sigma): E(|g|)\in
L^0(\mathcal{A})\}$. Throughout this paper we take $u$ and $w$ in
$\mathcal{D}(E)$. \\

Every operator $T$ on a Hilbert space $\mathcal{H}$ can be
decomposed into $T = U|T|$ with a partial isometry $U$, where $|T|
= (T^*T)^{\frac{1}{2}}$ . $U$ is determined uniquely by the kernel
condition $\mathcal{N}(U) = \mathcal{N}(|T|)$.  Then this
decomposition is called the polar decomposition. The Aluthge
transformation of $T$ is the operator $\hat{T}$ given by
$\hat{T}=|T|^{\frac{1}{2}}U|T|^{\frac{1}{2}}$. \\
The plan for the reminder of this paper is to present
characterizations of weighted conditional expectation type
operators in some various normality classes. Here is a brief
review of what constitutes membership for an operator $T$ on a
Hilbert space in each
classes:\\
(i) $T$ is normal if $T^{\ast}T=TT^{\ast}$.\\
%(ii) $T$ is quasinormal if $T(T^{\ast}T)=(T^{\ast}T)T$.\\
(ii) $T$ is hyponormal if $T^{\ast}T\geq TT^{\ast}$.\\
(iii) For $0<p<\infty$, $T$ is p-hyponormal if $(T^{\ast}T)^p\geq (TT^{\ast})^p$.\\
(iv) $T$ is $\infty$-hyponormal if it is p-hyponormal for all p.\\
%(Vi) $T$ is quasihyponormal if $T(T^{\ast}T)\geq(T^{\ast}T)T$.\\
%(Vii) For $M>0$, $T$ is $M$-quasihyponormal if
%$M^2T^{\ast^2}T^2\geq (T^{\ast}T)^2$.\\
(v) $T$ is weakly hyponormal if
$|\hat{T}|\geq|T|\geq|\hat{T}^{\ast}|$.\\
%(iX) $T$ belongs to $\ast$-$A$-class if
%$|T^2|\geq|T^{\ast}|^2$.\\
%(X) $T$ belongs to $\ast$-$A$-class if
%$T^{\ast}|T^2|T\geq T^{\ast}|T^{\ast}|^2T$.\\
%(Xi) The operator $T$ belongs to class $A$ if $|T|^2\leq|T^2|$.\\
(vi) $T$ is normaloid if $\|T\|^n=\|T^n\|$ for all
$n\in\mathbb{N}$.

%There are several basic relationships between these classes. The
%ones of concern in this note are as follows:

%\begin{center}
%normal $\Rightarrow$ quasinormal $\Rightarrow$ $\infty$-hyponormal
%$\Rightarrow$ $p$-hyponormal $\Rightarrow$ weakly hyponormal
%$\Rightarrow$ $p$-quasihyponormal $\Rightarrow$ normaloid;
%\end{center}
%\begin{center}
%$p$-hyponormal $\Rightarrow$ $p$-quasihyponormal $\Rightarrow$
%normaloid.
%\end{center}

\section{ \sc\bf Some classes of weighted conditional expectation type operators}

In the first we reminisce some theorems that we have proved in
\cite{ej}.

\vspace*{0.3cm} {\bf Theorem 3.1.} The operator $T=M_wEM_u$ is
bounded on $L^{2}(\Sigma)$ if and only if
$(E|w|^{2})^{\frac{1}{2}}(E|u|^{2})^{\frac{1}{2}} \in
L^{\infty}(\mathcal{A})$, and in this case its norm is given by
$\|T\|=\|(E(|w|^{2}))^{\frac{1}{2}}(E(|u|^{2}))^{\frac{1}{2}}\|_{\infty}$.\\

\vspace*{0.3cm} {\bf Lemma 3.2.} Let $T=M_wEM_u$ be a bounded
operator on $L^{2}(\Sigma)$ and let $p\in (0,\infty)$. Then
$$(T^{\ast}T)^{p}=M_{\bar{u}(E(|u|^{2}))^{p-1}\chi_{S}(E(|w|^{2}))^{p}}EM_{u}$$
and
$$(TT^{\ast})^{p}=M_{w(E(|w|^{2}))^{p-1}\chi_{G}(E(|u|^{2}))^{p}}EM_{\bar{w}},$$
where $S=S(E(|u|^2))$ and $G=S(E(|w|^2))$.

\vspace*{0.3cm} {\bf Theorem 3.3.} The unique polar decomposition
of bounded operator $T=M_wEM_u$ is $U|T|$, where

$$|T|(f)=\left(\frac{E(|w|^{2})}{E(|u|^{2})}\right)^{\frac{1}{2}}\chi_{S}\bar{u}E(uf)$$
and
 $$U(f)=\left(\frac{\chi_{S\cap
 G}}{E(|w|^{2})E(|u|^{2})}\right)^{\frac{1}{2}}wE(uf),$$
for all $f\in L^{2}(\Sigma)$.

\vspace*{0.3cm} {\bf Theorem 3.4.} The Aluthge transformation of
$T=M_wEM_u$ is
$$\widehat{T}(f)=\frac{\chi_{S}E(uw)}{E(|u|^{2})}\bar{u}E(uf), \ \ \ \ \ \  \ \ \ \  \ \  \ f\in
L^{2}(\Sigma).$$\\

 From now on, we consider the operators $M_wEM_u$ and $EM_u$ are
 bounded operators on $L^{2}(\Sigma)$.
 In the sequel some necessary and sufficient conditions for
 normality, hyponormality, $p$-hyponormality, . . . will be presented.\\

\vspace*{0.3cm} {\bf Theorem 3.5.} \vspace*{0.3cm} Let
$T=M_wEM_u$, then

(a) If
$(E(|u|^2))^{\frac{1}{2}}\bar{w}=u(E(|w|^2))^{\frac{1}{2}},$ then
$T$ is normal.

(b) If $T$ is normal, then $|E(u)|^2E(|w|^2)=|E(w)|^2E(|u|^2)$.\\

{\bf Proof.} (a) Applying lemma 3.2 we have
$$T^{\ast}T-TT^{\ast}=M_{\bar{u}E(|w|^{2})}EM_{u}-M_{wE(|u|^{2})}EM_{\bar{w}}.$$
So for every $f\in L^{2}(\Sigma)$,
$$\langle T^{\ast}T-TT^{\ast}(f),f\rangle=$$
$$\int_{X}E(|w|^2)E(uf)\bar{uf}-E(|u|^2)E(\bar{w}f)w\bar{f}d\mu$$
$$=\int_{X}|E(u(E(|w|^2))^{\frac{1}{2}}f)|^2-|E((E(|u|^2))^{\frac{1}{2}}\bar{w}f)|^2d\mu.$$
This implies that if
$$(E(|u|^2))^{\frac{1}{2}}\bar{w}=u(E(|w|^2))^{\frac{1}{2}},$$
then for all $f\in L^2(\Sigma)$, $\langle
T^{\ast}T-TT^{\ast}(f),f\rangle=0$, thus $T^{\ast}T=TT^{\ast}$.\\

(b) Suppose that $T$ is normal. For all $f\in L^2(\Sigma)$ we have
$$\int_{X}|E(u(E(|w|^2))^{\frac{1}{2}}f)|^2-|E((E(|u|^2))^{\frac{1}{2}}\bar{w}f)|^2d\mu=0.$$
Let $A\in \mathcal{A}$, with $0<\mu(A)<\infty$. By replacing $f$
to $\chi_{A}$, we have

$$\int_{A}|E(u(E(|w|^2))^{\frac{1}{2}})|^2-|E((E(|u|^2))^{\frac{1}{2}}\bar{w})|^2d\mu=0$$
and so
$$\int_{A}|E(u)|^2E(|w|^2)-|E(w)|^2E(|u|^2)d\mu=0.$$
Since $A\in \mathcal{A}$ is arbitrary, then
$|E(u)|^2E(|w|^2)-|E(w)|^2E(|u|^2)$.$\Box$

\vspace*{0.3cm} {\bf Corollary 3.6.} The operator $EM_u$ is normal
if and only if $u\in L^{\infty}(\mathcal{A})$.

 \vspace*{0.3cm} {\bf Theorem 3.7.}
  Let $T=M_wEM_u$ and let $p\in (0,\infty)$.\\
(a) The followings are equivalent.\\

(i) $T$ is hyponormal.\\

(ii) $T$ is $p$-hyponormal.\\

(iii) $T$ is $\infty$-hyponormal.\\

(b) If $u(E(|w|^2))^{\frac{1}{2}}-(E(|u|^2))^{\frac{1}{2}}\bar{w}\geq0,$ then $T$ is hyponormal.\\

(c) If $T$ is hyponormal, then
$|E(u)|^2E(|w|^2)-|E(w)|^2E(|u|^2)\geq0$.\\

\vspace*{0.3cm} {\bf Proof.} (a) Applying Lemma 3.2 we obtain that
$(T^{\ast}T)^{p}\geq(TT^{\ast})^{p}$ if and only if

$$M_{\chi_{S\cap
 G}(E(|u|^{2}))^{p-1}(E(|w|^{2}))^{p-1}}(M_{\bar{u}E(|w|^{2})}EM_{u}-M_{wE(|u|^{2})}EM_{\bar{w}})\geq0.$$

This inequality holds if and only if
$$T^{\ast}T-TT^{\ast}=M_{\bar{u}E(|w|^{2})}EM_{u}-M_{wE(|u|^{2})}EM_{\bar{w}}\geq 0,$$
where we have used the fact that $T_{1}T_{2}\geq0$ if
$T_{1}\geq0$, $T_{2}\geq0$ and $T_{1}T_{2}=T_{2}T_{1}$ for all
$T_{i}\in \mathcal{B}(\mathcal{H})$, the set of all bounded linear
operators on Hilbert space $\mathcal{H}$. Since $0<p<\infty$ is
arbitrary, then (i), (ii) and (iii) are equivalent.\\

(b) By lemma 3.2 we have
$$T^{\ast}T-TT^{\ast}=M_{\bar{u}E(|w|^{2})}EM_{u}-M_{wE(|u|^{2})}EM_{\bar{w}}.$$
So for every $f\in L^{2}(\Sigma)$,
$$\langle T^{\ast}T-TT^{\ast}(f),f\rangle=$$
$$=\int_{X}|E(u(E(|w|^2))^{\frac{1}{2}}f)|^2-|E((E(|u|^2))^{\frac{1}{2}}\bar{w}f)|^2d\mu.$$

This implies that, if $$
 u(E(|w|^2))^{\frac{1}{2}}-(E(|u|^2))^{\frac{1}{2}}\bar{w}\geq0,$$ then $T$ is hyponormal.

 (c) Let $T$ be hyponormal. For all $f\in L^2(\Sigma)$
we have
$$\int_{X}|E(u(E(|w|^2))^{\frac{1}{2}}f)|^2-|E((E(|u|^2))^{\frac{1}{2}}\bar{w}f)|^2d\mu\geq0.$$
Let $A\in \mathcal{A}$, with $0<\mu(A)<\infty$. By replacing $f$
to $\chi_{A}$, we have

$$\int_{A}|E(u(E(|w|^2))^{\frac{1}{2}})|^2-|E((E(|u|^2))^{\frac{1}{2}}\bar{w})|^2d\mu\geq0$$
and so
$$\int_{A}|E(u)|^2E(|w|^2)-|E(w)|^2E(|u|^2)d\mu\geq0.$$
Since $A\in \mathcal{A}$ is arbitrary, then
$|E(u)|^2E(|w|^2)\geq|E(w)|^2E(|u|^2)$. $\Box$

  \vspace*{0.3cm} {\bf Corollary 3.8.}
  Let $T=EM_u$, and $p\in (0,\infty)$. Then the followings are equivalent.\\

(i) $T$ is normal.\\

(ii) $T$ is hyponormal.\\

(iii) $T$ is $p$-hyponormal.\\

(iV) $T$ is $\infty$-hyponormal.\\

(V) $u\in L^{\infty}(\mathcal{A})$.

\vspace*{0.3cm} {\bf Theorem 3.9.} Let $T=M_wEM_u$, then\\

 (a) If $|E(uw)|^{2}\geq
E(|u|^{2})E(|w|^{2})$, then $T$ is
$p$-quasihyponormal.\\

(b) If $T$ is $p$-quasihyponormal, then
 $|E(uw)|^{2}\geq E(|u|^{2})E(|w|^{2})$ on $\sigma(E(u))\cap G$.\\

 (c) If $S(w)=S(u)=X$, then  $T$ is
 $p$-quasihyponormal if and only if $|E(uw)|^{2}\geq E(|u|^{2})E(|w|^{2})$.

 \vspace*{0.3cm} {\bf Proof.}
(a) By Lemma 3.2, it is easy to check that
$$\ \ \ \ \ \ \ \ \ \ \ T^{\ast}(T^{\ast}T)^{p}T=M_{\bar{u}(E(|u|^{2}))^{p-1}\chi_{S}(E(|w|^{2}))^{p}|E(uw)|^{2}}EM_{u};$$
$$T^{\ast}(TT^{\ast})^{p}T=M_{\bar{u}(E(|w|^{2}))^{p+1}(E(|u|^{2}))^{p}}EM_{u} .$$

It follows that $T^{\ast}(T^{\ast}T)^{p}T\geq
T^{\ast}(TT^{\ast})^{p}T$ if
$$M_{(E(|u|^{2}))^{p-1}\chi_{S}(E(|w|^{2}))^{p}}M_{(|E(uw)|^{2}-E(|w|^{2})E(|u|^{2}))}M_{\bar{u}}EM_{u}\geq0.$$

By the same argument in Theorem 3.7, this inequality holds if
$M_{(|E(uw)|^{2}-E(|w|^{2})E(|u|^{2}))}\geq 0$; i.e.
$|E(uw)|^{2}-E(|w|^{2})E(|u|^{2}) \geq0$.

\vspace*{0.3cm} (b) Suppose that  $T$ is $p$-quasihyponormal. Then
for all $f\in L^{2}(\mathcal{A})$, we have

$$\langle T^{\ast}(T^{\ast}T)^{p}T-T^{\ast}(TT^{\ast})^{p}Tf,f\rangle
$$$$
=\int_{X}(E(|u|^{2}))^{p-1}\chi_{S}(E(|w|^{2}))^{p}(|E(uw)|^{2}-E(|w|^{2})E(|u|^{2}))|E(u)|^{2}|f|^{2}d\mu\geq0.$$
Thus
$$(E(|u|^{2}))^{p-1}\chi_{S}(E(|w|^{2}))^{p}(|E(uw)|^{2}-E(|w|^{2})E(|u|^{2}))|E(u)|^{2}\geq0 ,$$
and hence we obtain $|E(uw)|^{2}\geq E(|w|^{2})E(|u|^{2})$ on $\sigma(E(u))\cap G$.

\vspace*{0.3cm} (c) It follows from (a) and (b).$\Box$

\vspace*{0.3cm} {\bf Corollary 3.10.}
  Let $T=EM_u$, $S(u)=X$ and $p\in (0,\infty)$. Then the following cases are equivalent.\\

(i) $T$ is hyponormal.\\

(ii) $T$ is $p$-hyponormal.\\

(iii) $T$ is $\infty$-hyponormal.\\

(iV) $T$ is $p$-quasihyponormal.\\

(V) $u\in L^{\infty}(\mathcal{A})$.

\vspace*{0.3cm} {\bf Example 3.11.} Let $X=[0,1]\times [0,1]$,
$d\mu=dxdy$, $\Sigma$  the  Lebesgue subsets of $X$ and let
$\mathcal{A}=\{A\times [0,1]: A \ \mbox{is a Lebesgue set in} \
[0,1]\}$. Then, for each $f$ in $L^2(\Sigma)$, $(Ef)(x,
y)=\int_0^1f(x,t)dt$, which is independent of the second
coordinate. This example is due to A. Lambert and B. Weinstock
\cite{la}. Now, if we take $u(x,y)=y^{\frac{x}{8}}$ and $w(x,
y)=\sqrt{(4+x)y}$, then $E(|u|^2)(x,y)=\frac{4}{4+x}$ and
$E(|w|^2)(x,y)=\frac{4+x}{2}$. So, $E(|u|^2)(x,y)E(|w|^2)(x,y)=2$
and $|E(uw)|^2(x,y)=64\frac{4+x}{(x+12)^2}$. Direct computations
shows that $E(|u|^2)(x,y)E(|w|^2)(x,y)\leq|E(uw)|^2(x,y)$. Thus,
by Theorem 3.9 the weighted conditional type operator $M_wEM_u$ is
$p$-quasihyponormal.

\vspace*{0.3cm} {\bf Theorem 3.12.} Let $T=M_wEM_u$, then\\

 (a) If $|E(uw)|=E(|u|^2)(E(|w|^2))^{\frac{1}{2}}$ on $S=S(E(|u|^2))$,
then $T$ is weakly hyponormal.\\

(b) If $T$ is weakly hyponormal, then
$|E(uw)|=E(|u|^2)(E(|w|^2))^{\frac{1}{2}}$ on $S(E(u))$.

 \vspace*{0.3cm} {\bf Proof.} (a)
For every $f\in L^2(\Sigma)$ by Theorem 3.3 and 3.4 we have

$$|\widehat{T}|(f)=|(\widehat{T})^{\ast}|(f)=|E(uw)|\chi_{S}(E(|u|^2))^{\frac{-3}{2}}\bar{u}E(uf),$$
where $S=S(E(|u|^2))$.

So, $T$ is weakly hyponormal if and only if $|T|=|\widehat{T}|$.
For every $f\in L^2(\Sigma)$,
$$\langle|T|(f)-|\widehat{T}|(f),f\rangle=
\int_{X}\left(\frac{E(|w|^{2})}{E(|u|^{2})}\right)^{\frac{1}{2}}\chi_{S}\bar{uf}E(uf)-|E(uw)|\chi_{S}(E(|u|^2))^{\frac{-3}{2}}\bar{uf}E(uf)d\mu$$
$$\int_{X}\left(\frac{E(|w|^{2})}{E(|u|^{2})}\right)^{\frac{1}{2}}\chi_{S}|E(uf)|^2-|E(uw)|\chi_{S}(E(|u|^2))^{\frac{-3}{2}}|E(uf)|^2d\mu,$$
this implies that if $|E(uw)|=E(|u|^2)(E(|w|^2))^{\frac{1}{2}}$ on
$S$, then $|T|=|\widehat{T}|$.\\

(b) If $|T|=|\widehat{T}|$, then for all $f\in L^2(\Sigma)$ we
have

$$\int_{X}\left(\frac{E(|w|^{2})}{E(|u|^{2})}\right)^{\frac{1}{2}}\chi_{S}|E(uf)|^2-|E(uw)|\chi_{S}(E(|u|^2))^{\frac{-3}{2}}|E(uf)|^2d\mu=0.$$

Let $A\in \mathcal{A}$, with $0<\mu(A)<\infty$. By replacing $f$
to $\chi_{A}$, we have

$$\int_{A}\left(\frac{E(|w|^{2})}{E(|u|^{2})}\right)^{\frac{1}{2}}\chi_{S}|E(u)|^2-|E(uw)|\chi_{S}(E(|u|^2))^{\frac{-3}{2}}|E(u)|^2d\mu=0.$$

Since $A\in \mathcal{A}$ is arbitrary, then
$$\left(\frac{E(|w|^{2})}{E(|u|^{2})}\right)^{\frac{1}{2}}\chi_{S}|E(u)|^2-|E(uw)|\chi_{S}(E(|u|^2))^{\frac{-3}{2}}|E(u)|^2=0.$$

Hence $|E(uw)|=E(|u|^2)(E(|w|^2))^{\frac{1}{2}}$ on
$S(E(u))$.$\Box$

\vspace*{0.3cm} {\bf Corollary 3.13.}\\

(a) If $T=EM_u$ and $S(E(u))=S(E(|u|^2))$, then $T$ is weakly
hyponormal if and only if $|E(u)|=E(|u|^2)$ on $S(E(u))$.\\

(b) If $T=M_wE$, then $T$ is weakly hyponormal if and only if
$w\in L^{\infty}(\mathcal{A})$.\\

\section{ \sc\bf Some Applications}

From now on, we shall denote by $\sigma(T)$, $\sigma_{p}(T)$,
$\sigma_{jp}(T)$, $\sigma_{ap}(T)$, $r(T)$ the spectrum of $T$,
the point spectrum of $T$, the approximate point spectrum, the
joint point spectrum of $T$, the spectral radius of $T$,
respectively. The spectrum of an operator $T$ is the set
$$\sigma(T)=\{\lambda\in \mathbb{C}:T-\lambda I \ \ \  {\text{is \  not \
invertible}}\}.$$ A complex number $\lambda\in \mathbb{C}$ is said
to be in the point spectrum $\sigma_{p}(T)$ of the operator $T$,
if there is a unit vector $x$ satisfying $(T-\lambda)x=0$. If in
addition, $(T^{\ast}-\bar{\lambda})x=0$, then $\lambda$ is said to
be in the joint point spectrum $\sigma_{jp}(T)$ of $T$. The approximate point spectrum of $T$ is the set of those $\lambda$ such that $T-\lambda I$
is not bounded below. Also, the spectral radius of $T$ is defined by $r(T)=\sup\{|\lambda|: \lambda\in \sigma(T)\}$.\\

If $A, B\in \mathcal{B}(\mathcal{H})$, then it is well known that
$$\sigma(AB)\setminus\{0\}=\sigma(BA)\setminus\{0\},$$
$$\sigma_{p}(AB)\setminus\{0\}=\sigma_{p}(BA)\setminus\{0\},$$
$$\sigma_{ap}(AB)\setminus\{0\}=\sigma_{ap}(BA)\setminus\{0\},$$
$$\sigma_{jp}(AB)\setminus\{0\}=\sigma_{jp}(BA)\setminus\{0\}.$$\\

 J. Herron
showed that if $EM_u:L^2(\Sigma)\rightarrow L^2(\Sigma)$, then
$\sigma(EM_u)=ess \ range(E(u))\cup\{0\}$, \cite{her}. So we
conclude that
$$\sigma(M_wEM_u)\setminus \{0\}=ess \
range(E(uw))\setminus \{0\}.$$ Let $A_{\lambda}=\{x\in
X:E(u)(x)=\lambda\}$, for $0\neq\lambda\in \mathbb{C}$. Suppose
that $\mu(A_{\lambda})>0$. Since $\mathcal{A}$ is $\sigma$-finite,
there exists an $\mathcal{A}$-measurable subset $B$ of
$A_{\lambda}$ such that $0<\mu(B)<\infty$, and $f=\chi_{B}\in
L^p(\mathcal{A})\subseteq L^p(\Sigma)$. Now
$$EM_u(f)-\lambda f=E(u)\chi_{B}-\lambda \chi_{B}=0.$$ This
implies that $\lambda\in P_{\sigma}(EM_u)$.\\
If there exists $f\in L^p(\Sigma)$ such that $f\chi_{C}\neq 0$
$\mu$-a.e, for $C\in \Sigma$ of positive measure and
$E(uf)=\lambda f$ for $0\neq \lambda \in \mathbb{C}$, then
$f=\frac{E(uf)}{\lambda}$, which means that $f$ is
$\mathcal{A}$-measurable. Therefore $E(uf)=E(u)f=\lambda f$ and
$(E(u)-\lambda)f=0$. This implies that $C\subseteq A_{\lambda}$
and so $\mu(A_{\lambda})>0$. Hence
$$\sigma_p(EM_u)=\{\lambda\in\mathbb{C}:\mu(A_{\lambda})>0\}.$$
Thus
$$\sigma_p(M_wEM_u)\setminus \{0\}=\{\lambda\in\mathbb{C}:\mu(A_{\lambda,w})>0\}\setminus \{0\},$$
where $A_{\lambda,w}=\{x\in X:E(uw)(x)=\lambda\}$.\\

For each natural number $n$, we define
$$\triangle_n(T)=\widehat{\triangle_{n-1}{T}} \ \ \ \ \
\triangle_1(T)=\triangle(T)=\widehat{T}.$$

We call $\triangle_n(T)$ the $n$-th Aluthge transformation of $T$.
It is proved that $r(T)=\lim_{n\rightarrow
\infty}\|\triangle_n(T)\|$ in \cite{ty}.\\

\vspace*{0.3cm} {\bf Theorem 4.1.}  Let $T=M_wEM_u$. Then\\

(a) $\widehat{T}$ is normaloid.\\

(b) $T$ is normaloid if and only if
$$\|E(uw)\|_{\infty}=\|(E(|u|^2))^{\frac{1}{2}}(E(|w|^2))^{\frac{1}{2}}\|_{\infty}.$$

\vspace*{0.3cm} {\bf Proof.} (a) By Theorem 3.1 we have
$\|\widehat{T}\|=\|E(uw)\|_{\infty}$. By Theorem 3.4 we conclude
that for every natural number $n$ we have
$\triangle_n(T)=\triangle(T)=\widehat{T}$. Hence
$r(\widehat{T})=r(T)=\|\widehat{T}\|=\|E(uw)\|_{\infty}$. So
$\widehat{T}$ is normaloid.\\

(b) By conditional type Holder inequality, boundedness of $T$ and
Theorem 3.1 we have
$$r(T)=\|E(uw)\|_{\infty}\leq\|(E(|u|^2))^{\frac{1}{2}}(E(|w|^2))^{\frac{1}{2}}\|_{\infty}=\|T\|.$$

Hence $T$ is normaloid if and only if
$$\|E(uw)\|_{\infty}=\|(E(|u|^2))^{\frac{1}{2}}(E(|w|^2))^{\frac{1}{2}}\|_{\infty}.$$

 $\Box$
\vspace*{0.3cm} {\bf Theorem 4.2.}  Let $T=M_wEM_u$ be weakly
hyponormal with $ker T\subset ker T^{\ast}$, then
$T=\widehat{T}$.\\

\vspace*{0.3cm} {\bf Proof.} Direct computations show that
$\widehat{T}$ is normal. So, by Theorem 2.6 of \cite{aw} we have
$T=\widehat{T}$.

\vspace*{0.3cm} {\bf Theorem 4.3.}  If $|E(uw)|=E(|u|^2)(E(|w|^2))^{\frac{1}{2}}$ on $S=S(E(|u|^2))$, then\\

(a)  $\sigma_{jp}(M_wEM_u)\setminus \{0\}=\sigma_p(M_wEM_u)\setminus \{0\}=\{\lambda\in\mathbb{C}:\mu(A_{\lambda,w})>0\}\setminus \{0\}$.\\

(b) $\sigma(M_wEM_u)\setminus \{0\}=\{\lambda:\bar{\lambda}\in
\sigma_{ap}(M_{\bar{u}}EM_{\bar{w}})\}\setminus\{0\}$.\\

(c) If $\lambda$ is an isolated point in $\sigma(M_wEM_u)$, then
$\lambda\in \sigma_p(M_wEM_u)$.\\

\vspace*{0.3cm} {\bf Proof.} By using Theorem 3.12 and Theorems
3.2 and 3.7 of
\cite{aw} we conclude (a), (b) and (c).\\

\vspace*{0.3cm} {\bf Corollary 4.4.} If $|E(uw)|=E(|u|^2)(E(|w|^2))^{\frac{1}{2}}$ on $S=S(E(|u|^2))$, then\\

$$\sigma_{ap}(M_wEM_u)\setminus
\{0\}=\sigma(M_wEM_u)\setminus \{0\}=ess \ range(E(uw))\setminus
\{0\}.$$

%If $T$ is paranormal, then
%$$E(|w|^2)|E(u)|^2\leq |E(uw)|(E(|w|^2))^{\frac{1}{2}}|E(u)|.$$

\end{document}